\documentclass{amsart}
\usepackage{amsmath, amssymb, latexsym}
\title{Bounding Betti numbers of bipartite graph ideals}
\author{Michael Goff}

\newtheorem{theorem}{Theorem}[section]
\newtheorem{proposition}[theorem]{Proposition}

\newtheorem{problem}[theorem]{Problem}

\newcommand{\K}{\Gamma}
\newcommand{\field}{{\bf k}}
\newcommand{\codim}{\mbox{\upshape codim}\,}

\newcommand{\mindeg}{\mbox{\upshape mindeg}\,}

\def\proof{\smallskip\noindent {\it Proof: \ }}
\def\proofof#1{\smallskip\noindent {\it Proof of #1: \ }}
\def\endproof{\hfill$\square$\medskip}

\begin{document}

\begin{abstract}
We prove a conjectured lower bound of Nagel and Reiner on Betti numbers of edge ideals of bipartite graphs.
\end{abstract}

\date{July 9, 2008}

\maketitle

\section{Introduction and preliminaries}
Finding explicit minimal free resolutions for classes of graded ideals, or at least bounding their Betti numbers, is one of the central problems in combinatorial algebra.  In general, the problem is hard and far from being solved, even in the cases of monomial ideals or quadratic monomial ideals (for some results and conjectures, see e.g. \cite{LinSyz}, \cite{HVT}, and the survey paper \cite{EdgeIdealSurvey}).  In this paper we prove a conjecture raised by Nagel and Reiner \cite{NR}, establishing a lower bound on the Betti numbers of certain quadratic ideals.

We start by reviewing necessary background and introducing notation.  Throughout this paper $\field$ is an arbitrary field, and $S$ is the polynomial ring over $\field$ in variables $X \sqcup Y$, where $X = \{x_1,\ldots,x_n\}$ and $Y=\{y_1,\ldots,y_m\}$.  We consider ideals generated by some monomials of the form $x_iy_j$.  Define a $\mathbb{Z}^n$-grading on $S$ as follows.  Let $\mathbb{Z}^n$ be generated by the standard basis $e_1,\ldots,e_n$, and set $\deg x_i = e_i$ for $1 \leq i \leq n$.  Also set $\deg y_i=0$ for $1 \leq i \leq m$.

For a $\mathbb{Z}^n$-graded ideal $I \subset S$, we consider the minimal free $\mathbb{Z}^n$-graded resolution:
\[ 0 \rightarrow
\bigoplus_{a \in \mathbb{Z}^n} S(-a)^{\beta_{l,a}} \rightarrow \ldots \rightarrow \bigoplus_{a \in \mathbb{Z}^n} 
S(-a)^{\beta_{0,a}} \rightarrow I \rightarrow 0. \]
In the above expression, $S(-a)$ denotes $S$ with grading shifted by $a$, and $l$ denotes the length of the resolution. In particular, $l \geq \codim(S/I)$.
It follows from, for instance, the Taylor resolution that if $I$ is a squarefree monomial ideal, then $\beta_{i,a}=0$ unless $a$ is a $\{0,1\}$-vector.  Hence the nonzero Betti numbers of such an ideal can be indexed by subsets of $X$.  For $X' \subseteq X$, we define $\beta_{i,X',\bullet}(I) = \beta_{i,a}(I)$ for $a = \sum_{x_i \in X'}e_i$.  We may also consider the more common $\mathbb{Z}^{n+m}$-grading on $S$ by giving $y_j$ degree $e_{n+j}$.  Then the $\mathbb{Z}^n$-graded Betti numbers of $I$ relate to the $\mathbb{Z}^{n+m}$-graded Betti numbers  by

\begin{equation}
\label{zn}
\beta_{i,X',\bullet}(I) = \sum_{Y' \subseteq Y}\beta_{i,X' \sqcup Y'}(I).
\end{equation}

In Section \ref{LowerBound}, we prove Conjecture 1.2 of \cite{NR}, establishing a lower bound on $\beta_{i,X',\bullet}(I)$ in the case that $I$ is generated by some monomials of the form $x_iy_j$.  Associated with $I$ is a bipartite graph $G(X \sqcup Y,E)$ with vertex set $X \sqcup Y$ and an edge $\{x_i,y_j\} \in E$ whenever $x_iy_j \in I$.  We say that $I$ is the \textit{edge ideal} of $G$.  Edge ideals were first introduced in \cite{EdgeIdeals}; results related to edge ideals can be found in \cite{Splittable}, \cite{Jacques}, \cite{Forests}, \cite{CharInd}, and \cite{EdgeIdeals}.  For each vertex $v \in G$, the set of vertices that share an edge with $v$ is called the \textit{neighborhood} of $v$ and is denoted $N(v)$, while the \textit{degree} of $v$ is $\deg v = \deg_G v := |N(v)|$.  

For each bipartite graph $G$ on $X \sqcup Y$, we associate a bipartite graph $H$ on $X \sqcup Y$ with edge set
$$E(H) = \{\{x_i,y_j\}: 1 \leq i \leq n,\;\; 1 \leq j \leq \deg_G x_i\}. $$
We may think of $H$ as a ``shifted'' version of $G$.  A bipartite graph constructed in this manner is known as a \textit{Ferrers graph}.  Let $J$ be the edge ideal of $H$; $J$ is known as a \textit{Ferrers ideal}.  For more on Ferrers ideals, see \cite{Ferrers} and \cite{Ferrers2}.  The following is Conjecture 1.2 of \cite{NR}.

\begin{theorem}
\label{bipartite}
For all $X' \subset X$, $\beta_{i,X',\bullet}(I) \geq \beta_{i,X',\bullet}(J).$
\end{theorem}

Our proof relies heavily on techniques relating to simplicial complexes.  A \textit{simplicial complex} $\K$ with the vertex set $V = X \sqcup Y$ is a collection of subsets of $2^{V}$ called \textit{faces} such that if $F \in \K$ and $G \subseteq F$, then $G \in \K$.  With every simplicial complex $\K$ we associate its \textit{Stanley-Reisner ideal} $I_{\K} \subset S$ generated by non-faces of $\K$: $I_\K := (\prod_{v \in L}v: L \subseteq V, L \not\in \K)$ (see \cite{St96}).  Likewise, given a squarefree monomial ideal $I \subset S$, we denote by $\Delta(I)$ the simplicial complex $\Delta$ on $X \sqcup Y$ such that $I_\Delta = I$.  If $W \subset V$, then the \textit{induced subcomplex} of $\K$ on $W$, denoted $\K[W]$ has vertex set $W$ and faces $\{F \in \K: F \subseteq W\}$.  If $v \in V$, then we abbreviate $\K[V-\{v\}]$ by $\K-v$.  Let $\tilde{\beta}_p(\K) := \dim_\field(\tilde{H}_p(\K))$ be the dimension of the $p$-th reduced simplicial homology of $\K$ with coefficients in $\field$.  We make frequent use of Hochster's formula (see \cite[Theorem II.4.8]{St96}), which states that for $W \subset V$, $$\beta_{i,W}(I_\K) = \tilde{\beta}_{i-|W|-2}(\K[W]).$$

\section{Lower bound on bipartite graph ideals}
\label{LowerBound}
In this section we prove the main result.  Let $G$ be a graph on $X \sqcup Y$, all of whose edges are of the form $\{x_i,y_j\}$, and let $I$ be the edge ideal of $G$.  Let $J$ be the Ferrers ideal associated with $I$.  The Betti numbers of Ferrers ideals can be calculated explicitly.  For $X' \subseteq X$, let $\mindeg(X') = \mindeg_G(X')$ denote the minimum degree of a vertex in $X'$ in $G$.  

\begin{proposition}
\label{FerrersBN}
\cite[Proposition 2.18]{NR}.  Let $J$ be the edge ideal of a Ferrers graph $H$ on vertex set $X \sqcup Y$.  Then for all $X' \subseteq X$ and $i$, $$\beta_{i,X',\bullet}(J) = {\mindeg_H(X') \choose i-|X'|+2}.$$  
\end{proposition}

\proofof{Theorem \ref{bipartite}} For a given $X' \subseteq X$, we may restrict our attention to the induced subgraph $G[X' \sqcup Y]$, and therefore we assume without loss of generality that $X' = X$.  By Proposition \ref{FerrersBN}, $\beta_{i,X,\bullet}(J) = {\mindeg(X) \choose i-|X|+2}$.  Let $\K := \Delta(I)$.  By (\ref{zn}) and Hochster's formula, we also have that $$\beta_{i,X,\bullet}(I) = \sum_{Y' \subseteq Y}\beta_{i,X \sqcup Y'}(I) = \sum_{j=0}^{|X|+|Y|-i-2}\sum_{|Y'| = j+i-|X|+2}\tilde{\beta}_{j}(\K[X \cup Y']).$$

We assume without loss of generality that $N(x_1)$ does not properly contain $N(x_i)$ for $1 \leq i \leq n$.  This occurs, for instance, if $x_1$ has minimal degree among the vertices in $X$.  It suffices to show that $$\sum_{j=0}^{|X|+|Y|-i-2}\sum_{|Y'| = j+i-|X|+2}\tilde{\beta}_{j}(\K[X \cup Y']) \geq {\deg(x_1) \choose i-|X|+2}.$$  We do so by showing that for every $Y_1' \subset N(x_1)$, there exists $Y' \subseteq Y$ and $j \geq 0$ such that $Y' \cap N(x_1) = Y_1'$, $|Y'|=|Y_1'|+j$, and $\tilde{\beta}_j(\K[X \cup Y']) \geq 1$.  If this claim holds, then by taking all $Y_1'$ with $|Y_1'|=i-|X|+2$, it follows that $\beta_{i,|X|,\bullet}(I) \geq {\deg(x_1) \choose i-|X|+2}$.

Define $X_1 := \{x \in X: N(x) = N(x_1)\}$.  If $x \in X-X_1$, then there exists $y \in N(x)-N(x_1)$, since by our hypothesis $N(x) \not\subset N(x_1)$.  Let $\{v_1, \ldots, v_r\} \subset Y-N(x_1)$ be a set of minimal size such that for each $x \in X-X_1$, there exists some $1 \leq i \leq r$ with $v_i \in N(x)$.  We prove the claim by induction on $r$.  In the case $r=0$, $N(x_i)=N(x_1)$ for all $i$, and so $\K[X \cup Y_1'] = \K[X_1 \cup Y_1']$ is the disjoint union of simplices on $X_1$ and $Y_1'$, and the claim holds with $j=0$.

Now consider $r \geq 1$, and let $X' = N(v_r)$.  On the induced subgraph $G[(X-X') \cup Y'_1 \cup \{v_1,\ldots,v_{r-1}\}]$, $N(x_1)$ does not properly contain $N(x_i)$ for any $x_i \in X-X'$, so for this graph the claim holds by the inductive hypothesis.  Hence by possibly rearranging the $v_i$, we can assume that $\tilde{H}_{k-1}(\K[(X-X') \cup Y'_1 \cup \{v_1,\ldots,v_{k-1}]\}) \neq 0$ for some $1 \leq k \leq r$.  Then we consider two cases.

Case 1: $\tilde{H}_{k-1}(\K[X \cup Y'_1 \cup \{v_1,\ldots,v_{k-1}\}]) \neq 0$.  Then $Y' = Y_1' \cup \{v_1, \ldots, v_{k-1}\}$ satisfies the claim.

Case 2: $\tilde{H}_{k-1}(\K[X \cup Y'_1 \cup \{v_1,\ldots,v_{k-1}\}]) = 0$.  Note that $$\tilde{H}_{k-1}(\K[(X-X') \cup Y'_1 \cup \{v_1,\ldots,v_{k-1},v_r\}]) = 0$$ since this complex is a cone over $\K[(X-X') \cup Y'_1 \cup \{v_1,\ldots,v_{k-1}\}]$ with apex $v_r$.  Also, since for all $x \in X'$, $\{x,v_r\}$ is not an edge in $\K$, it follows that 
\begin{eqnarray*}
& & \K[X \cup Y'_1 \cup \{v_1,\ldots,v_{k-1},v_r\}] = \\
& & \K[X \cup Y'_1 \cup \{v_1,\ldots,v_{k-1}\}] \cup \K[(X-X') \cup Y'_1 \cup \{v_1,\ldots,v_{k-1},v_r\}].\\
\end{eqnarray*}
Take $X^* := X-X'$.  The portion of the Mayer-Vietoris sequence on simplicial homology
\begin{eqnarray*}
& & \tilde{H}_k(\K[X \cup Y'_1 \cup \{v_1, \ldots,v_{k-1}, v_r\}]) \rightarrow \tilde{H}_{k-1}(\K[X^* \cup Y'_1 \cup \{v_1, \ldots, v_{k-1}\}]) \rightarrow \\
& & \tilde{H}_{k-1}(\K[X \cup Y'_1 \cup \{v_1, \ldots, v_{k-1}\}]) \oplus \tilde{H}_{k-1}(\K[X^* \cup Y'_1 \cup \{v_1, \ldots,v_{k-1}, v_r\}])=0 
\end{eqnarray*}
implies that $\tilde{H}_k(\K[X \cup Y'_1 \cup \{v_1, \ldots,v_{k-1}, v_r\}]) \neq 0$.  The result follows by taking $Y' = Y_1' \cup \{v_1, \ldots,v_{k-1}, v_r\}$.
\endproof

Nagel and Reiner give a full characterization of when equality occurs for all $X' \subseteq X$.  We say that $G$ is \textit{nearly row-nested} if whenever $|N(x_1)| < |N(x_2)|$, $N(x_1) \subset N(x_2)$, and $|\cap_{|N(x_i)|=c} N(x_i)| \geq c-1$ for all $c$.

\begin{theorem}
\label{mainEQ}
\cite[Proposition 4.18]{NR} For all $X' \subset X$, $\beta_{i,X',\bullet}(I) = \beta_{i,X',\bullet}(J)$ if and only if $G$ is nearly row-nested.
\end{theorem}

\section{Remarks and conclusions}
Nagel and Reiner also propose a \textit{colex lower bound} for classes of monomial ideals.  The \textit{colex order} on subsets of size $d$ of $\mathbb{N}$ is a total ordering such that $(a_1,\ldots,a_d) <_{colex} (b_1,\ldots,b_d)$ if and only if for some $1 \leq k \leq d$, $a_k < b_k$ and $a_i = b_i$ for all $k+1 \leq i \leq d$.  An initial segment $K$ in the colex order is a \textit{colexsegment}, and the ideal $(x_{i_1}\ldots x_{i_d}: \{i_1,\ldots,i_d\} \in K)$ is a \textit{colexsegment-generated} ideal.  For each squarefree monomial ideal $I$ generated in a constant degree $d$, let $J$ be the unique degree $d$ colexsegment-generated ideal with the same number of minimal generators as $I$.  We say that $I$ satisfies the colex lower bound if for all $j$, $\beta_j(I) \geq \beta_j(J)$.  Problem 1.1 of \cite{NR} is the following.

\begin{problem}
Which monomial ideals in constant degree $d$ satisfy the colex lower bound?
\end{problem}

Theorem \ref{bipartite} proves the colex lower bound for edge ideals of bipartite graphs.

\begin{theorem}
Let $G$ be a bipartite graph.  Then the edge ideal of $G$ satisfies the colex lower bound.
\end{theorem}
\proof Let $I$ be the edge ideal of $G$ and $J$ be the associated Ferrers ideal.  Nagel and Reiner \cite[Proposition 4.2]{NR} prove that $J$ satisfies the colex lower bound.  By ignoring the $\mathbb{Z}^n$-grading, it follows from Theorem \ref{bipartite} that for all $j$, $\beta_j(I) \geq \beta_j(J)$.  We conclude that $I$ satisfies the colex lower bound.
\endproof

\end{document}